\documentclass[12pt]{article}
\usepackage{amsmath,amsthm,amssymb,latexsym,color}
\usepackage{bbm}

\usepackage{graphicx}

\usepackage{hyperref}

\date{}

\def\kkk{\null\hfill $\Box $ \\}

\newtheorem{theorem}{Theorem}[section]

\newcommand{\E}{\mathbb E\,}
\newcommand{\R}{\mathbb{R}}

\newcommand{\B}{\mathbb{B}}

\newcommand{\conv}{\mathop{\mathrm{conv}}\nolimits}

\newcommand{\ind}{\mathbbm{1}}

\newcommand{\bx}{\mathbf{x}}

\newcommand{\dd}{{\rm d}}

\newcommand{\El}{\mathcal{E}}

\newcommand\qand{\quad \mbox{ and } \quad }

\title{Random section and random simplex
inequality.}
\author{Alexander E. Litvak\thanks{
The work of this author was
supported by  Ministry of Science and Higher Education
of the Russian Federation, agreement No. 075-15-2019-1619.
}
and Dmitry Zaporozhets\thanks{
The work of this author  was supported by RFBR and DFG according to the research project No. 20-51-12004 and by the Foundation for the
Advancement of Theoretical Physics and Mathematics ``BASIS''.
}}

\newcommand\address{\noindent\leavevmode

\medskip
\noindent
Alexander E. Litvak\\
Dept.~of Math.~and Stat.~Sciences,\\
University of Alberta, \\
Edmonton, AB, Canada, T6G 2G1.\\
\texttt{\small
e-mail:  aelitvak@gmail.com}\\

\medskip

\noindent
Dmitry Zaporozhets\\
St.~Petersburg Department of\\
Steklov Institute of Mathematics\\
St.~Petersburg, Russia\\
\texttt{\small
e-mail:   zap1979@gmail.com}\\
}

\begin{document}

\maketitle

\begin{abstract}
Consider some convex body $K\subset\mathbb R^d$. Let  $X_1,\dots, X_k$, where $k\leq d$, be random points independently and uniformly chosen in $K$, and let $\xi_k$ be a uniformly distributed  random linear $k$-plane. We show that for $p\geq -d+k+1$,
\[
\mathbb E\,|K\cap\xi_k|^{d+p}\leq c_{d,k,p}\cdot|K|^k\,\,\mathbb E\,|\mathrm{conv}(0,X_1,\dots,X_k)|^p,
\]
where $|\cdot|$ and $\mathrm{conv}$ denote the volume of correspondent dimension and the convex hull. The constant $c_{d,k,p}$ is  such that for $k>1$ the equality holds if and only if $K$ is an ellipsoid centered at the origin, and for $k=1$ the inequality turns to equality.

If $p=0$, then the inequality reduces to the Busemann intersection inequality, and if $k=d$ -- to the Busemann random simplex inequality.

We also present an affine version of this inequality which similarly generalizes the Schneider inequality and the Blaschke-Gr\"omer inequality.

\end{abstract}



\maketitle


\section{Introduction}
\label{s:intro}

\subsection{Busemann intersection inequality}

For $d\in\mathbb N$ and $k\in\{1,\dots,d\}$, the linear Grassmannian of $k$-dimensional linear subspaces of $\R^d$ is denoted by $G_{d,k}$  and is equipped with a unique rotation invariant probabilistic Haar measure $\nu_{d,k}$. By $|\cdot|$ we denote the $d$-dimensional volume. Given $k\leq d$,
slightly abusing notation, considering the sets intersected with $k$-dimensional
affine subspaces or the convex hulls of $k+1$ points, we denote the $k$-dimensional volume
by $|\cdot|$ as well.

The seminal Busemann intersection inequality states that for any  convex compact set  $K\subset\R^d$ with non-empty interior (i.e., a convex body),
\begin{align}\label{1552}
   \int\limits_{G_{d,k}}|K\cap L|^{d}\,\nu_{d,k}(\dd L)\leq \frac{\kappa_k^d}{\kappa_d^k}\, |K|^{k},
\end{align}
where $\kappa_k$ denotes the volume of the $k$-dimensional unit ball.

Originally Busemann~\cite{hB53}  proved this inequality for  $k=d-1$ and  later it was generalized in \cite{BS60,eG92}  for all  $k=1,\dots,d-1$. If $k>1$,	 then the equality holds if and only if $K$ is an ellipsoid centered at the origin, and in this case
\eqref{1552} turns to the classical Fustenberg--Tzkoni formula \cite{FT71}.

Using the polar coordinates, it is easy to see that for $k=1$ the inequality  turns to the equality:
\begin{align}\label{0949}
   \int\limits_{G_{d,1}}|K\cap L|^{d}\,\nu_{d,1}(\dd L)=\frac{2^d}{\kappa_d}\, |K|.
\end{align}

Moreover, this equation can be  generalized to other moments as follows:
\begin{align}\label{1613}
    \int_{G_{d,1}}|&K\cap L|^{d+p}\,\nu_{d,k}(\dd L)
=\frac{(d+p)2^{d+p}}{d\kappa_d}
\int_{K^{k}}{|\bx|^p\, \dd\bx},\quad p\geq -d+k+1.
\end{align}

\noindent
{\bf Problem 1.1} {\it Is it possible for all  $k=1,\dots,d-1$ to obtain a generalization of~\eqref{1552} of the form
\begin{align*}
   \int\limits_{G_{d,k}}|K\cap L|^{d+p}\,\nu_{d,k}(\dd L)\leq\dots
\end{align*}
which turns to~\eqref{1613} for $k=1$?}

Now let us consider an affine version of~\eqref{1552}. To this end, denote by $A_{d,k}$  the affine Grassmannian of $k$-dimensional affine subspaces of $\R^d$ equipped with a unique    measure $\nu_{d,k}$ invariant with respect to the rigid motions in $\R^d$ and normalized by
\[
\mu_{d,k}\left(\left\{E\in A_{d,k}:\,E\cap \B^d\ne\emptyset\right\}\right)=\kappa_{d-k}.
\]
Schneider~\cite{rS85} showed that
\begin{align}\label{1600}
   \int\limits_{A_{d,k}}|K\cap E|^{d+1}\,\mu_{d,k}(\dd E)\leq \frac{\kappa_k^{d+1}\kappa_{d(k+1)}}{\kappa_d^{k+1}\kappa_{k(d+1)}}\, |K|^{k+1},
\end{align}
and for $k>1$ the equality holds if and only if $K$ is an ellipsoid.

As above, for $k=1$ the inequality turns to equality, although it is not as trivial as in the linear case, see
\cite{mC85} for $d=2$ and \cite{hH52} for any $d$. As in the linear case, this equality can be generalized to
other moments. It was done independently in \cite[Eq.~(21)]{gC67} and \cite[Eq.~(34)]{jK69}: for $p\geq -d+k+1$,
\begin{align}\label{1244}
    \int\limits_{A_{d,1}}|K\cap E|^{p+d+1}\,\mu_{d,1}(dE)=\frac{(d+p)\,(d+p+1)}{2d\kappa_d}\int\limits_{K^2}{|\bx_0-\bx_1|^p\,\dd\bx_0\dd\bx_1}.
\end{align}

\noindent
{\bf Problem 1.2} {\it As in the linear case, it is natural to ask: is it possible for all  $k=1,\dots,d-1$ to obtain a generalization of~\eqref{1600} of the form
\begin{align*}
    \int\limits_{A_{d,k}}|K\cap E|^{d+p+1}\,\mu_{d,k}(\dd E)\leq \dots,
\end{align*}
which turns to~\eqref{1244} for $k=1$?}

\medskip

To conclude this section, let us note that Gardner~\cite{rG07} generalized~\eqref{1552} and~\eqref{1600} to bounded Borel sets and characterized the equality cases. Recently Dann, Paouris, and Pivovarov~\cite{DPP16} extended \eqref{1552},~\eqref{1600} to bounded integrable functions.

\subsection{Busemann random simplex inequality}
Another  group of inequalities deals with the volume of the random simplex in a body. The classical Busemann random simplex inequality states that
\begin{align*}
	|K|^{d+1}\leq(d+1)!\frac{\kappa_d^{d+1}}{2\kappa_{d+1}^{d-1}}\int\limits_{K^{d}}{|\conv(0,\bx_1,\ldots,\bx_d)|\,\dd\bx_1\ldots \dd\bx_d}.
\end{align*}
This inequality can be generalized (see, e.g.,~\cite[Theorem~8.6.1.]{SW08}) as follows: for every $p\geq 1$,
\begin{align}\label{1027}
 |K|^{p+d} \leq (d!)^{p} \, \frac{\kappa_d^{p+d}}{\kappa_{d+p}^{d}}\, b_{d+p,d}\,
 \int\limits_{K^{d}} {|\conv(0,\bx_1,\ldots,\bx_d)|^p\, \dd\bx_1\ldots \dd\bx_d},
\end{align}
where for  a non-integer $p>0$ we denote
\begin{equation}\label{1409}
\kappa_p:=\frac{\pi^{p/2}}{\Gamma\left(\frac p2+1\right)}\quad \qand \quad b_{q,k}:={q\choose k}\frac{\kappa_{q-k+1}\cdots\kappa_{q}}{\kappa_1\cdots\kappa_k}.
\end{equation}
The equality holds if and only if $K$ is a centered ellipsoid.

\noindent
{\bf Problem 1.3} {\it Is it possible  to obtain a generalization of~\eqref{1027} for a random simplex of arbitrary dimension $k$?}

The affine counterpart of~\eqref{1027} is known as the Blaschke-Gr\"omer inequality~\cite{hG73}:
for every $p\geq 1$,
\begin{align}\label{1045}
    |K|^{p+d+1}\leq(d!)^{p}b_{d+p,d}\frac{\kappa_d^{p+d+1}}{\kappa_{d+p}^{d+1}}\,\frac{\kappa_{(d+1)(d+p)}}{\kappa_{d(d+p+1)}}
\int\limits_{K^{d+1}}{|\conv(\bx_0,\ldots,\bx_d)|^p\,\dd\bx_0\ldots \dd\bx_d}.
\end{align}
The equality holds if and only if $K$ is an ellipsoid.

\
\noindent
{\bf Problem 1.4} {\it  As in the linear case, it is natural to ask: Is it possible  to obtain a generalization of~\eqref{1045} for a random simplex of arbitrary dimension $k$?}

\medskip

The aim of this  note is to derive  a general inequality which simultaneously solves Problems~1.1 and~1.3 thus generalizing both Busemann inequalities in one form.  

  positively solve  Problems 1.1 and 1.3  thus presenting a general form for  both  the Busemann intersection inequality and the Busemann random simplex inequality. We also present an affine version of the inequality which solves Problems~1.2 and~1.4 and thus implies  the Schneider inequality and the Blaschke-Gr\"omer inequality.

\section{Main results}
Our first theorem generalizes~\eqref{1552} and~\eqref{1027}.

\begin{theorem}\label{1114}
	For any convex body $K\subset\R^d$, $k\in\{0,1,\dots,d\}$, and any real number $p\geq -d+k+1$,
\begin{align}\label{1141}
    \int_{G_{d,k}}|&K\cap L|^{p+d}\,\nu_{d,k}(\dd L)&
\leq(k!)^{p}\, \, \frac{\kappa_k^{d+p}}{\kappa_{d+p}^{k}}\,\,
\frac{ b_{d+p,k}}{ b_{d,k}}
\int_{K^{k}}{|\conv(0,\bx_1,\ldots,\bx_k)|^p\, \dd\bx_1\ldots \dd\bx_k}.
\end{align}
  For $k>1$  the equality holds if and only if $K$ is a non-degenerate  ellipsoid centered at the origin.
\end{theorem}

\noindent
{\bf Remarks.}
\begin{enumerate}
\item
Applying (\ref{1141}) with $p=0$ we obtain \eqref{1552}, while
applying it with $k=d$ we obtain \eqref{1027}.
\item
It was shown in~\cite[Theorem~1.6]{GGZ19} that if $K$ is a non-degenerate
ellipsoid centered at the origin,
then one has the equality in \eqref{1141}.
\item
In the probabilistic language it may be formulated as
\[
\E|K\cap\xi_k|^{p+d}\leq (k!)^{p}\, \, \frac{\kappa_k^{d+p}}{\kappa_{d+p}^{k}}\,\,
\frac{ b_{d+p,k}}{ b_{d,k}} \,\,|K|^k\,\,\E|\conv(0,X_1,\dots,X_k)|^p,
\]
where $X_1,\dots,X_k$ are independently and uniformly distributed points in $K$ and $\xi_k$ is a random linear $k$-plane uniformly distributed in $G_{d,k}$.
\end{enumerate}
\medskip

Our second theorem generalizes~\eqref{1600} and~\eqref{1045}.
\begin{theorem}\label{1146}
	For any convex body $K\subset\R^d$, $k\in\{0,1,\dots,d\}$, and any real number $p\geq -d+k+1$,
	\begin{align}\label{1152}
        \int_{A_{d,k}}&|K\cap E|^{p+d+1}\,\mu_{d,k}(dE)
         \leq C(k, p, d)
        \int_{K^{k+1}}{|\conv(\bx_0,\ldots,\bx_k)|^p\,\dd\bx_0\ldots \dd\bx_k},
    \end{align}
    where
$$
     C(k, p, d) = (k!)^{p}\, \frac{\kappa_k^{p+d+1}}{\kappa_{d+p}^{k+1}}\,
         \frac{\kappa_{(k+1)(d+p)}}{\kappa_{k(d+p)+k}}\,\frac{b_{d+p,k}}{b_{d,k}}.
$$
     For $k>1$  the equality holds if and only if $K$ is a non-degenerate ellipsoid.
\end{theorem}

{\bf Remarks.}
\begin{enumerate}
\item
Applying (\ref{1152}) with $p=0$ we obtain \eqref{1600}, while
applying it with $k=d$ we obtain \eqref{1045}.
\item
It was shown in~\cite[Theorem~1.4]{GGZ19} that if $K$ is a non-degenerate ellipsoid,
then one has the equality in~\eqref{1152}.
\item
In probabilistic language~\eqref{1152} may be formulated as
\begin{align*}
    \E|K\cap\eta_k|^{p+d+1}\leq  C'(k, p, d) \,\,
        \frac{|K|^{k+1}}{V_{d-k}(K)}\,\,\E|\conv(X_0,X_1,\dots,X_k)|^p,
\end{align*}
where
\begin{align*}
    C'(k, p, d) = \frac{d!\,(k!)^{p-1}}{(d-k)!}\,\frac{\kappa_d}{\kappa_{d-k}} \frac{\kappa_k^{p+d}}{\kappa_{d+p}^{k+1}}\,
         \frac{\kappa_{(k+1)(d+p)}}{\kappa_{k(d+p)+k}}\,\frac{b_{d+p,k}}{b_{d,k}},
\end{align*}
 $X_0,X_1,\dots,X_k$  are independently and uniformly distributed points in $K$, $\eta_k$ is uniformly distributed among all affine $k$-planes intersected $K$, and $V_{d-k}$ is the $(d-k)$-th intrinsic volume of $K$ defined by the Crofton formula~\cite[Theorem~5.1.1]{SW08} as the normalized measure of all affine $k$-planes intersected $K$:
\begin{align*}
    V_{d-k}(K):=\binom{d}{k}\frac{\kappa_d}{\kappa_k\kappa_{d-k}}\mu_{d,k}\left(\left\{E\in A_{d,k}:\,E\cap K\ne\emptyset\right\}\right).
\end{align*}
\end{enumerate}

\section{Proofs}\label{037}

\subsection{Blaschke--Petkantschin formula}\label{2257}

Recall that $b_{d,k}$ is defined by \eqref{1409}.
Given points $\bx_0, \bx_1,\ldots,\bx_k\in \R^d$ we denote
$$
  V_k = V(\bx_0, \bx_1,\ldots,\bx_k) := |\conv(\bx_0, \bx_1,\ldots,\bx_k)|
$$
and
$$
 V_{0,k} = V(\bx_1,\ldots,\bx_k):= |\conv(0,\bx_1,\ldots,\bx_k)|.
$$
In our further calculations we will need to integrate some non-negative measurable
 function $h$ of $k$-tuples of points in $\R^d$. To this end, we first integrate
 over the $k$-tuples of points in a fixed $k$-dimensional linear subspace $L$ and then we integrate over $G_{d,k}$. The corresponding transformation formula is known
 as the linear Blaschke--Petkantschin formula (see \cite[Theorem 7.2.1]{SW08}):
\begin{align}\label{eq3}
	\int\limits_{(\R^d)^k}&{h\, \dd\bx_1\ldots \dd\bx_k}=
     (k!)^{d-k}b_{d,k}
     \int\limits_{G_{d,k}}\int\limits_{L^k}h\,
       V_{0,k}^{d-k}\,\lambda_L(\dd\bx_1)\ldots\lambda_L(\dd\bx_k)\,\nu_{d,k}(\dd L),
\end{align}
where   $h=h(\bx_1,\ldots,\bx_k)$.
The following is an  affine counterpart of  (\ref{eq3}):
\begin{align}\label{eq3_2}
	\int\limits_{(\R^d)^{k+1}}&{h\,\dd\bx_0\ldots \dd\bx_k}=
   (k!)^{d-k}b_{d,k} \int\limits_{A_{d,k}}\int\limits_{E^{k+1}}h\,
   V_{ k }^{d-k}\,\lambda_E(\dd\bx_0)\ldots\lambda_E(\dd\bx_k)\,
   \mu_{d,k}(\dd E),
\end{align}
where  $h=h(\bx_0, \bx_1,\ldots,\bx_k)$
(see~\cite[Theorem 7.2.7]{SW08}).

\subsection{Proof of Theorem~\ref{1114}}\label{sec012}

Let

$$
	J:=\int\limits_{K^{k}}{V_{0,k}^p\,\dd\bx_1\ldots \dd\bx_k}
	=\int\limits_{(\R^d)^{k}}{V_{0,k}^{p}\,\prod\limits_{i=1}^k\ind_{K}(\bx_i)\,\dd\bx_1\ldots \dd\bx_k}.
$$
Applying the linear Blaschke--Petkantschin formula \eqref{eq3} with the function
\[
h(\bx_1,\ldots,\bx_k):=V_{0,k}^{p}\,\prod\limits_{i=1}^k\ind_{K}(\bx_i),
\]
we get
\begin{align}\label{91}
   J&=(k!)^{d-k}b_{d,k}\,\int\limits_{G_{d,k}}\int\limits_{L^{k}}
   V_{0,k}^{p+d-k}\,\prod\limits_{i=1}^k\ind_{K}(\bx_i)\,
   \lambda_L(\dd\bx_1)\ldots\lambda_L(\dd\bx_k)\,\nu_{d,k}(\dd L)\notag\\
	&=(k!)^{d-k}b_{d,k}\,\int\limits_{G_{d,k}}\int\limits_{(K\cap L)^{k}}
    V_{0,k}^{p+d-k}\,\lambda_L(\dd\bx_1)\ldots\lambda_L(\dd\bx_k)\,\nu_{d,k}(\dd L).
\end{align}
Fix $L\in G_{d,k}$. Applying \eqref{1027} with $p+d-k$ and $k$ instead of $p$ and $d$, we obtain
\begin{align}\label{1310}
	(k!)^{p+d-k}\, \frac{\kappa_{k}^{d+p}}{\kappa_{d+p}^{k}}\, b_{d+p, k}\,
     \int\limits_{(K\cap L)^{k}}V_{0,k}^{p+d-k}\,
      \lambda_L(\dd\bx_1)\ldots\lambda_L(\dd\bx_k) \geq |K\cap L|^{p+d},
\end{align}
which together with~\eqref{91} implies \eqref{1141}.

Finally we consider the equality case. As was mentioned above, the equality holds for ellipsoids, see~\cite[Theorem~1.6]{GGZ19}.
Conversely, suppose that~\eqref{1141} turns to equality. Then it follows from~\eqref{91} that~\eqref{1310} turns to equality
for almost all $L\in G_{d,k}$ which, in fact, means that it is true for \emph{all} $L\in G_{d,k}$. Indeed, if for some
$L\in G_{d,k}$ we had a strict inequality in~\eqref{1310}, then the same would be true for some neighborhood of $L$ which
would contradict to the fact that~\eqref{1310} turns to equality for almost all $L\in G_{d,k}$.
Thus, according to the equality case in (\ref{1027}),  $K\cap L$ is a centered
ellipsoid for all $L\in G_{d,k}$.
Now it remains to apply
the following lemma from~\cite[(16.12)]{hB55}: if for some fixed $k>1$ for any $E\in A_{d,k}$ passing through some fixed point from the
interior of $K$ the intersection $K\cap E$ happens to be a $k$-dimensional ellipsoid, then $K$ is an ellipsoid itself.
\kkk

\subsection{Proof of Theorem~\ref{1146}}
The proof is  similar to the previous one.  Let
$$
	J:=\int\limits_{K^{k+1}}{V_k^p\,\dd\bx_0\ldots \dd\bx_k}
	=\int\limits_{(\R^d)^{k+1}}{V_k^{p}\,\prod\limits_{i=0}^k\ind_{\El}(\bx_i)\,\dd\bx_0\ldots \dd\bx_k}.
$$
Applying  the affine Blaschke--Petkantschin formula \eqref{eq3_2} with the function
\[
h(\bx_0,\ldots,\bx_k):=|\conv(\bx_0,\ldots,\bx_k)|^{p}\,\prod\limits_{i=0}^k\ind_{\El}(\bx_i),
\]
we get
\begin{align}\label{1931}
  J&=(k!)^{d-k}b_{d,k}\,\int\limits_{A_{d,k}}
  \int\limits_{E^{k+1}}V_k^{p+d-k}\,
   \prod\limits_{i=0}^k\ind_{K}(\bx_i)\,\lambda_E(\dd\bx_0)\ldots\lambda_E(\dd\bx_k)\,\mu_{d,k}(\dd E)
   \notag\\
	&=(k!)^{d-k}b_{d,k}\,\int\limits_{A_{d,k}}\int\limits_{(K\cap E)^{k+1}}V_k^{p+d-k}\,
   \lambda_E(\dd\bx_0)\ldots\lambda_E(\dd\bx_k)\,\mu_{d,k}(\dd E).
\end{align}
Fix $E\in A_{d,k}$. Applying \eqref{1045} with $p+d-k$ and $k$ instead of $p$ and $d$, we obtain
\begin{align*}
   (k!)^{d-k+p}b_{d+p,k}\, \frac{\kappa_k^{p+d+1}}{\kappa_{d+p}^{k+1}}\,
   \frac{\kappa_{(k+1)(d+p)}}{\kappa_{k(d+p+1)}}\,  \int\limits_{(K\cap E)^{k+1}}V_k^{p+d-k}\,\lambda_E(\dd\bx_0)\ldots\lambda_E(\dd\bx_k)
    \\ \geq |K\cap E|^{d+p+1}
\end{align*}
which together with~\eqref{1931} implies \eqref{1152}. The equality case is treated the same way as in the linear case.
\kkk

\bibliography{bib}
\bibliographystyle{plain}

\address

\end{document}